\documentclass[journal,onecolumn]{IEEEtran}

\usepackage{cite}
\usepackage{indentfirst}
\usepackage{amsmath,amssymb,amsfonts}
\usepackage{upgreek}
\renewcommand{\footnoterule}{%
  \kern -3pt
  \hrule width \textwidth height 0.4pt
  \kern 2pt
}
\newcommand\mycom[2]{\genfrac{}{}{0pt}{}{#1}{#2}}

\makeatletter
\newcommand*{\rom}[1]{\expandafter\@slowromancap\romannumeral #1@}
\makeatother

\ifCLASSINFOpdf
\else
\fi
\begin{document}
\title{Relative Performance of Fisher Information in Interval Estimation}

\author{Sihang Jiang \\ {Department of Engineering Systems and Environment} \\
{University of Virginia}\\
Charlottesville, United States \\ {sj5yq@virginia.edu}}% <-this % stops a space
\maketitle

\begin{abstract}
Maximum likelihood estimates and corresponding confidence regions of the estimates are commonly used in statistical inference. In practice, people often construct approximate confidence regions with the Fisher information at given sample data based on the asymptotic normal distribution of the MLE (maximum likelihood estimate). Two common Fisher information matrices (FIMs, for multivariate parameters) are the observed FIM (the Hessian matrix of negative log-likelihood function) and the expected FIM (the expectation of the observed FIM). In this article, we prove that under certain conditions and with an MSE (mean-squared error) criterion, approximate confidence interval of each element of the MLE with the expected FIM is at least as accurate as that with the observed FIM.
\end{abstract}

\begin{IEEEkeywords}
MLE, observed FIM, expected FIM, interval estimation
\end{IEEEkeywords}

% \fi
\IEEEpeerreviewmaketitle

\section{Introduction}

Maximum likelihood estimates (MLE) and corresponding confidence regions are among the most popular methods for parameter estimation. The covariance matrix of the MLE is an important element while constructing the confidence region. Under reasonable conditions (Section 13.3 in Spall \cite{spall_2003}) the MLE has the following asymptotic normality:
\begin{equation}
\sqrt{n}(\boldsymbol{\hat{\uptheta}}_n-\boldsymbol{{\uptheta}}^*)\xrightarrow[\text{}]{\text{dist}} N(\boldsymbol{0},\boldsymbol{\bar{F}}(\boldsymbol{\uptheta}^*)^{-1})
\end{equation}
where $n$ is the sample size; $\boldsymbol{\uptheta}=[t_1,...,t_p]^T \in \mathbb{R}^p$ is the $p$-dimensional vector value parameter about which we are concerned; $\boldsymbol{\uptheta}^*=[t_1^*,...,t_p^*]^T$ denotes the true value that generates the data in nature; $\boldsymbol{\hat{\uptheta}}_n=[\hat{t}_1,...,\hat{t}_p]^T$ represents the MLE; and $\boldsymbol{\bar{F}}(\boldsymbol{\uptheta}^*)=\lim_{n\to\infty}\boldsymbol{F}_n(\boldsymbol{\uptheta}^*)/{n} $, with $\boldsymbol{F}_n(\boldsymbol{\uptheta}^*)$ denoting the Fisher information matrix (FIM) at $\boldsymbol{\uptheta}^*$. Due to the complexity of the setting, the asymptotic covariance matrix of the MLE $\boldsymbol{\bar{F}}(\boldsymbol{\uptheta}^*)^{-1}$, might be difficult or even impossible to obtain. In practice, if we have many collections of observations and accordingly get many simulations of the MLE, we can use the sample covariance matrix as an approximation. Also, the inverse of the Hessian matrix of the negative log-likelihood function $\boldsymbol{\bar{H}}(\boldsymbol{\uptheta}^*)^{-1}$, which is the inverse of observed FIM, is also considered as a possible approximation. Note that Guo and Spall \cite{guo_spall_2019} show that the inverse Hessian is superior to the sum of squared gradients of the log-likelihood as an "observed FIM" in the scalar case, and in this paper we show that under certain conditions the expected FIM gives a confidence level at least as accurate as the observed FIM\footnote{This paper is an expanded version of \cite{9400253}.}. \par
We now briefly review some research regarding which kind of FIM is a better choice. Efron and Hinkley \cite{10.2307/2335893} prove that for scalar parameter $\boldsymbol{\uptheta}^*$ if there exists an ancillary statistic $a$ (affects the precision of the MLE $\boldsymbol{\hat{\uptheta}}_n$), then the inverse of the observed Fisher information number (FIN) $\boldsymbol{\bar{H}}(\boldsymbol{\hat{\uptheta}}_n)^{-1}$ outperforms the expected FIN $\boldsymbol{\bar{F}}(\boldsymbol{\hat{\uptheta}}_n)^{-1}$ as an estimate for the conditional variance of the MLE $\textrm{var}(\boldsymbol{\hat{\uptheta}}_n|a)$. Lindsay and Li \cite{10.1214/aos/1069362393} compare the expected FIM and observed FIM as an approximation for the sample covariance matrix of MLE $\boldsymbol{\hat{\uptheta}}_n$. Specifically, they use following the MSE criterion: $$\min_{V} E[n(\boldsymbol{\hat{\uptheta}}_n-\boldsymbol{\uptheta}^*)(\boldsymbol{\hat{\uptheta}}_n-\boldsymbol{\uptheta}^*)^T-\boldsymbol{V}]^2$$ where $\boldsymbol{V}=\boldsymbol{\bar{F}}(\boldsymbol{\hat{\uptheta}}_n)^{-1}$ or $\boldsymbol{V}=\boldsymbol{\bar{H}}(\boldsymbol{\hat{\uptheta}}_n)^{-1}$. Moreover, they conclude that if an error term of magnitude $O(n^{-3/2})$ (with definition in Section \rom{2}.B) is ignored, the observed FIM gives the minimum MSE. Cao and Spall \cite{6315584} use the theoretical covariance matrix \cite{cao2013relative}, $n\textrm{cov}(\boldsymbol{\hat{\uptheta}}_n)$ as basis and show that the expected FIM is a better estimator asymptotically; i.e., $E[n\textrm{cov}(\boldsymbol{\hat{\uptheta}}_n)-\boldsymbol{\bar{F}}(\boldsymbol{\hat{\uptheta}}_n)^{-1}]^2 \leq E[n\textrm{cov}(\boldsymbol{\hat{\uptheta}}_n)-\boldsymbol{\bar{H}}(\boldsymbol{\hat{\uptheta}}_n)^{-1}]^2$. Yuan and Spall \cite{9147324} proved asymptotically, that the expected FIN outperforms the observed FIN in approximate confidence interval construction in an MSE sense in scalar cases. \par

\section{Mathematical Modeling}

\subsection{Problem Description}

Let $\boldsymbol{X}=[\boldsymbol{X}_1,\boldsymbol{X}_2,...,\boldsymbol{X}_n]$, $\boldsymbol{X}_i \in \mathbb{R}^q$, be a collection of independent but not necessarily identically distributed (i.n.i.d.) observations. Denote the probability density/mass functions (pdf/pmf) of $\boldsymbol{X}_i$ as 
$p_i(\boldsymbol{X}_i,\boldsymbol{\uptheta})$
for all $i=1,2,...,n$, where $\boldsymbol{\uptheta}$ is the vector parameter. The joint pdf/pmf is $$p(\boldsymbol{X},\boldsymbol{\uptheta})=\prod_{i=1}^{n}p_i(\boldsymbol{X}_i,\boldsymbol{\uptheta})$$ 
and denote the MLE of $\boldsymbol{\uptheta}$ as $\boldsymbol{\hat{\uptheta}}_n$ and the unknown true value as $\boldsymbol{\uptheta}^*$.\par

Specifically, for large $n$, $\boldsymbol{\hat{\uptheta}}_n$ obeys the following approximate normal distribution:
\begin{equation}
\sqrt{n}(\boldsymbol{\hat{\uptheta}}_n-\boldsymbol{{\uptheta}}^*) \sim N(\boldsymbol{0},\boldsymbol{V}_n)
\end{equation}
where $\boldsymbol{V}_n=n\textrm{cov}(\boldsymbol{\hat{\uptheta}}_n)$. In practice, $\boldsymbol{V}_n$ is approximated by either the expected or observed FIM.\par

Denote the negative log-likelihood function of $\boldsymbol{X_i}$  and $\boldsymbol{X}$ as 
$$l_i(\boldsymbol{X},\boldsymbol{\uptheta})\overset{\mathrm{\Delta}}{=}-\textrm{log}(p_i(\boldsymbol{X}_i,\boldsymbol{\uptheta}))$$
$$l(\boldsymbol{X},\boldsymbol{\uptheta})\overset{\mathrm{\Delta}}{=}\sum_{i=1}^{n}l_i(\boldsymbol{X},\boldsymbol{\uptheta})$$ and the observed FIM (denoted as $\boldsymbol{\bar{H}}_n$) and the expected FIM (denoted as $\boldsymbol{\bar{F}}_n$) are $$\boldsymbol{\bar{H}}_n(\boldsymbol{\uptheta},\boldsymbol{X})\overset{\mathrm{\Delta}}{=}n^{-1}{\partial^2 l(\boldsymbol{X},\boldsymbol{\uptheta})}/{\partial\boldsymbol{\uptheta}\partial\boldsymbol{\uptheta}^T}$$
$$\boldsymbol{\bar{F}}_n(\boldsymbol{\uptheta})\overset{\mathrm{\Delta}}{=}E[\boldsymbol{\bar{H}}_n(\boldsymbol{\uptheta},\boldsymbol{X})]$$ \par
Let $e_j$, $j=1,2,...,p$ be the unit $p$-dimensional vector with the $j^{th}$ component equal to 1 and the rest of the components equal to 0. Then, the $j^{th}$ component of the MLE $\boldsymbol{\hat{\uptheta}}_n$, $\hat{t}_{nj}$, satisfies $\hat{t}_{nj}=e_j^T\boldsymbol{\hat{\uptheta}}_n$ and $$\sqrt{n}(\hat{t}_{nj}-t_j^*) \sim N(0,\boldsymbol{V}_n(j,j))$$ where $\boldsymbol{V}_n(j,j)=e_j^T\boldsymbol{V}_ne_j$  is the $j^{th}$ diagonal component of $\boldsymbol{V}_n$.
Accordingly, the asymptotic $1-\upalpha_j$ confidence interval for $t_j^*$  (denoted as $CI_{jj}$) has the following expression: 
$$CI_{jj}=\left[\hat{t}_{nj}-z_{1-\upalpha_{j}/2}\sqrt{\frac{\boldsymbol{V}_n(j,j)}{n}},\hat{t}_{nj}+z_{1-\upalpha_{j}/2}\sqrt{\frac{\boldsymbol{V}_n(j,j)}{n}}\right]$$
where $\upalpha_j$ is the significance level with typical values like 0.01, 0.05, or 0.10; and $z_{1-\upalpha/2}$  is the corresponding standard normal deviate.
Substituting the ${jj}^{th}$ entry of the real covariance matrix with the corresponding entry of the inverse of the observed and the expected FIM (i.e. $\boldsymbol{\bar{H}}_n^{-1}(j,j)$ and $\boldsymbol{\bar{F}}_n^{-1}(j,j)$ for $\boldsymbol{V}_n(j,j)$), we derive the approximate $1-\upalpha_j$ confidence intervals (denoted as $CI_{H_{jj}}$ and $CI_{F_{jj}}$).\par

Based on the Bonferroni correction, let us choose $\upalpha_j$, $j=1,...,p$, such that $\upalpha=\sum_{j=1}^{p}\upalpha_j$, then the $p$-dimensional joint approximate $1-\upalpha$ confidence regions are as follows:
$$CI=CI_{11}\times...\times CI_{pp}$$
$$CI_H=CI_{H_{11}}\times...\times CI_{H_{pp}} $$
$$CI_F=CI_{F_{11}}\times...\times CI_{F_{pp}}$$
Moreover, for $t_j^*$, the approximate confidence levels with the observed and expected FIM are
$$
P\left(\frac{\sqrt{n}(\hat{t}_{nj}-t_j^*)}{\sqrt{\boldsymbol{V}_n(j,j)}}\in \left[-z_{1-\upalpha_j/2}\sqrt{\frac{\boldsymbol{\bar{H}}_n^{-1}(j,j)}{\boldsymbol{V}_n(j,j)}},z_{1-\upalpha_j/2}\sqrt{\frac{\boldsymbol{\bar{H}}_n^{-1}(j,j)}{\boldsymbol{V}_n(j,j)}}\right]\right)
$$
$$
P\left(\frac{\sqrt{n}(\hat{t}_{nj}-t_j^*)}{\sqrt{\boldsymbol{V}_n(j,j)}}\in \left[-z_{1-\upalpha_j/2}\sqrt{\frac{\boldsymbol{\bar{F}}_n^{-1}(j,j)}{\boldsymbol{V}_n(j,j)}},z_{1-\upalpha_j/2}\sqrt{\frac{\boldsymbol{\bar{F}}_n^{-1}(j,j)}{\boldsymbol{V}_n(j,j)}}\right]\right)
$$
so we define the confidence level function for the $j^{th}$ element of the parameter 
$$\uppi_j(x)\overset{\mathrm{\Delta}}{=}2\Phi \left(z_{1-\upalpha_j/2}\sqrt{\frac{x}{\boldsymbol{V}_n(j,j)}}\right)-1$$ 
then we have the approximate confidence level with the observed FIM: $$1-\upalpha_{H_{jj}}=P(t_j^* \in CI_{H_{jj}})=\uppi_j(\boldsymbol{\bar{H}}_n^{-1}(j,j))$$ and similarly, the approximate confidence level with the expected FIM is: $$1-\upalpha_{F_{jj}}=\uppi_j(\boldsymbol{\bar{F}}_n^{-1}(j,j))$$ and the asymptotic significance level is $$1-\upalpha_j=\uppi_j(\boldsymbol{V}_n(j,j))$$ \par
The confidence level represents the probability that the projection of the true parameter lies in the corresponding confidence interval. In subsequent parts of this article, we compare the mean squared errors (MSEs) of the confidence level of $t_j^*$ to compare the accuracy of the corresponding approximate confidence levels:

\begin{equation}
    \textrm{MSE}_{H_{jj}}=E[(1-\upalpha_j)-(1-\upalpha_{H_{jj}})]^2
\end{equation}
\begin{equation}
    \textrm{MSE}_{F_{jj}}=E[(1-\upalpha_j)-(1-\upalpha_{F_{jj}})]^2
\end{equation}

We show that under reasonable conditions, the MSE with the expected FIM ($\textrm{MSE}_{F_{jj}}$) is asymptotically smaller or equal to that with the observed FIM ($\textrm{MSE}_{H_{jj}}$). Further, the total MSE with the expected FIM ($\textrm{MSE}_{F_{jj}}$) is asymptotically smaller or equal to that with the observed FIM ($\textrm{MSE}_{H_{jj}}$). And these results imply that the approximate confidence interval with the expected FIM has no less accuracy compared to that with the observed FIM in the MSE sense for either any given element of the parameter or the entire parameter.

\subsection{Notation}
In this section, following the notations of Cao \cite{cao2013relative}, we define some random variables that are useful in the proof.\par
Denote the derivatives of the negative log-likelihood function as follows:
$$U_r\overset{\mathrm{\Delta}}{=}\frac{\partial l(\boldsymbol{\uptheta},x)}{\partial t_r}$$ $$U_{rs}\overset{\mathrm{\Delta}}{=}\frac{\partial^2 l(\boldsymbol{\uptheta},x)}{\partial t_r \partial t_s}$$ $$U_{rst}\overset{\mathrm{\Delta}}{=}\frac{\partial^3 l(\boldsymbol{\uptheta},x)}{\partial t_r \partial t_s \partial t_t}$$\par
Similarly, denote the null-cumulants as follows \cite{cao2013relative}:
$$\bar{\kappa}_r\overset{\mathrm{\Delta}}{=}n^{-1}E(U_r)=0$$ $$\bar{\upkappa}_{rs}\overset{\mathrm{\Delta}}{=}n^{-1}E(U_{rs})=\boldsymbol{\bar{F}}_n(\boldsymbol{\uptheta}^*)(r,s)$$ $$\bar{\upkappa}_{rst}\overset{\mathrm{\Delta}}{=}n^{-1}E(U_{rst})$$ $$\bar{\upkappa}_{r,s}\overset{\mathrm{\Delta}}{=}n^{-1}\textrm{cov}(U_r,U_s)$$
$$\bar{\upkappa}_{rs,t}\overset{\mathrm{\Delta}}{=}n^{-1}\textrm{cov}(U_{rs},U_t)$$
$$\bar{\upkappa}^{r,s}\overset{\mathrm{\Delta}}{=}\bar{\upkappa}_{rs}^{-1}$$ \par

Denote the standardized likelihood scores as follows: $$Z_r\overset{\mathrm{\Delta}}{=}n^{-{1}/{2}}U_r(\boldsymbol{\uptheta}^*)$$
$$Z_{rs}\overset{\mathrm{\Delta}}{=}n^{-{1}/{2}}[U_{rs}(\boldsymbol{\uptheta}^*)-n\bar{\upkappa}_{rs}]$$ $$Y_{rs}\overset{\mathrm{\Delta}}{=}Z_{rs}-\overline{\overline{\bar{\upkappa}_{rs,t}\bar{\upkappa}^{t,u}Z_u}}=\sum_{t=1}^p\sum_{u=1}^p\bar{\upkappa}_{rs,t}\bar{\upkappa}^{t,u}Z_u$$
where the double bar notation $\overline{\overline{(.)}}$ indicates a special summation operation. Specifically, for the argument under the double bar, summation is implied over any index repeated once as a superscript and once as a subscript. In addition, we frequently use the stochastic big-$O$ and little-$o$ terms: $O_d(n^{-r}), r \in \mathbb{R}$ represents a stochastic term that converges in distribution to a random variable when multiplied by $n^r$; $O_d^k(n^{-r})$  denotes the product of $k$ $O_d(n^{-r})$ terms; $O_d^2(n^{-r})$ denotes the square of an $O_d(n^{-r})$ term; $o_p(1)$ represents a stochastic term that converges in probability to zero; $o(1)$ represents a series of numbers that converge to zero; and $O_{+}(1)$ denotes a series of numbers that converge to a non-negative constant number. In addition, we introduce $\tilde{O_d}(n^{-r})$ to denote a summation of a finite number of $O_d(n^{-r})$ terms and $\tilde{O_d}^2(n^{-r})$ to denote a summation of a finite number of $O_d^2(n^{-r})$ terms.

\subsection{Conditions}

The main conclusion of this paper requires that some conditions be placed on the log-likelihood function and some sequences with respect to the log-likelihood function. These conditions include existence of partial derivatives, necessary interchanges of differentiation and integration, existence of limits and bounds, the weak law of large number and the dominated convergence theorem. A detailed discussion about the conditions may be found in Cao \cite[pp. 25\textendash28]{cao2013relative}. \par

We list an important condition which guarantees that the strict inequality in the main result holds: $\bar{\boldsymbol{F}}_n^{-1}(\hat{\boldsymbol{\uptheta}}_n)$ and  $\bar{\boldsymbol{H}}_n^{-1}(\hat{\boldsymbol{\uptheta}}_n,\boldsymbol{X})$ are not identical and we have
\begin{equation}
    \liminf\limits_{n\rightarrow \infty}n^{-1}\sum_{i=1}^n\textrm{var}[\overline{\overline{\bar{\upkappa}^{r,t}\bar{\upkappa}^{s,u}(U_{tu}^l-\bar{\upkappa}_{tu}^l-\bar{\upkappa}_{tu,v}\bar{\upkappa}^{v,w}U_w^l)}}]>0
\end{equation}

The summation term on the LHS of the inequality represents the main difference between $\bar{\boldsymbol{H}}_n^{-1}$  and $\bar{\boldsymbol{F}}_n^{-1}$. 

\subsection{Existing lemmas}

In this section, we present several lemmas from Cao \cite{cao2013relative} that illustrate the relationship between the variance of the MLE and the approximate variance with the two kinds of FIMs. For each cited lemma, we offer a brief explanation.\par
$\boldsymbol{Lemma \quad 1}$. For i.n.i.d sample data with conditions A1–A9 \cite{cao2013relative} holding, the estimation error for the components of $\hat{\boldsymbol{\uptheta}}_n$ have the following form:
$$\hat{t}_{nj}-t_j^*=-n^{-{1}/{2}}\overline{\overline{\bar{\upkappa}^{j,u}Z_u}}+\tilde{O_d}^2(n^{-1/2}), j=1,2,...p$$
Lemma 1 measures the difference between the MLE $\hat{\boldsymbol{\uptheta}}_n$ and the true parameter $\boldsymbol{\uptheta}^*$.\par
$\boldsymbol{Lemma \quad 2}$. For i.n.i.d sample data with conditions A1–A9 \cite{cao2013relative} holding, the variance of MLE   $\hat{\boldsymbol{\uptheta}}_n$ and the theoretical covariance matrix and the inverse of two FIMs have the following representation:
$$\textrm{var}(\hat{t}_{nj})=n^{-1}\bar{\upkappa}^{j,j}+o(n^{-1})$$
$$\bar{F}_n(\hat{\boldsymbol{\uptheta}}_n)^{-1}(j,j)=\bar{\upkappa}^{j,j}+B_{njj}+\tilde{O_d}^2(n^{-1/2})$$
$$\bar{H}_n(\hat{\boldsymbol{\uptheta}}_n)^{-1}(j,j)=\bar{\upkappa}^{j,j}+B_{njj}-A_{njj}+\tilde{O_d}^2(n^{-1/2})$$
where we have the expression $$A_{njj}=n^{-1/2}\overline{\overline{\bar{\upkappa}^{j,r}\bar{\upkappa}^{j,s}Y_{rs}}}$$
$$B_{njj}=n^{-1/2}\overline{\overline{\bar{\upkappa}^{j,t}\bar{\upkappa}^{j,u}\bar{\upkappa}^{v,w}(\bar{\upkappa}_{tuv}-\bar{\upkappa}_{tu,v})Z_w}}$$ and lemma 2 decomposes the variance of the $j^{th}$ element of MLE and the approximate variances with two kinds of FIMs and indicates their relationships to $\bar{\upkappa}^{j,j}$, the inverse of $\boldsymbol{\bar{F}}_n^{-1}(j,j)$.\par
$\boldsymbol{Lemma \quad 3}$. For i.n.i.d sample data with conditions A1–A9 \cite{cao2013relative} holding, and given $A_{njj}$  and $B_{njj}$ defined the same as in Lemma 2, we have:
$$E(A_{njj})=E(B_{njj})=E(A_{njj}B_{njj})=0$$
$$E[A_{njj}\tilde{O_d}^2(n^{-1/2})]=E[A_{njj}\tilde{O_d}^2(n^{-1/2})]=o(n^{-1})$$
$$E\{o(1)[-A_{njj}+\tilde{O_d}^2(n^{-1/2})]\}=o(n^{-1})$$
Lemma 3 presents properties of the difference terms $A_{njj}$  and $B_{njj}$. Specifically, $B_{njj}$ is the main difference between $\boldsymbol{\bar{F}}_n^{-1}(j,j)$ and  $\bar{\upkappa}^{j,j}$ while $B_{njj}-A_{njj}$  represents that between $\boldsymbol{\bar{H}}_n^{-1}(j,j)$ and $\bar{\upkappa}^{j,j}$.

\section{Theoretical analysis}
In this section, based on the existing lemmas, we offer relationships between relevant confidence intervals and corresponding confidence levels (a probability value). 

\subsection{Preliminary results}

$\boldsymbol{Lemma \quad 4}$. For i.n.i.d. sample data with conditions A1–A9 \cite{cao2013relative} holding, the convergence rates of the inverse of two FIMs to $\bar{\upkappa}^{j,j}$ are:
$$\boldsymbol{\bar{F}}_n(\hat{\boldsymbol{\uptheta}}_n)^{-1}(j,j)-\bar{\upkappa}^{j,j}=\tilde{O_d}(n^{-1/2})$$
$$\boldsymbol{\bar{H}}_n(\hat{\boldsymbol{\uptheta}}_n)^{-1}(j,j)-\bar{\upkappa}^{j,j}=\tilde{O_d}(n^{-1/2})$$

$Proof$: According to Condition A5, by the CLT for i.n.i.d. samples, $Z_r=n^{-1/2}U_r(\boldsymbol{\uptheta}^*)$ converges in distribution to a normal random variable. Thus $Z_r=\tilde{O_d}(1)$. Similarly,
$Z_{rs}=n^{-1/2}[U_{rs}(\boldsymbol{\uptheta}^*)-n\bar{\upkappa}_{rs}]=\tilde{O_d}(1)$.
It follows that $Y_{rs}$, the linear combination of  $Z_r$  and $Z_{rs}$, equals $\tilde{O_d}(1)$. Consequently, we have
$$A_{njj}=n^{-1/2}\overline{\overline{\bar{\upkappa}^{j,r}\bar{\upkappa}^{j,s}Y_{rs}}}=\tilde{O_d}(n^{-1/2})$$
$$B_{njj}=n^{-1/2}\overline{\overline{\bar{\upkappa}^{j,t}\bar{\upkappa}^{j,u}\bar{\upkappa}^{v,w}(\bar{\upkappa}_{tuv}-\bar{\upkappa}_{tu,v})Z_w}}=\tilde{O_d}(n^{-1/2})$$ and plugging expression of $A_{njj}$ and $B_{njj}$ into Lemma 2, we derive the result of Lemma 4.\par

$\boldsymbol{Lemma \quad 5}$. For i.n.i.d sample data with conditions A1–A9 \cite{cao2013relative} holding, the relevant confidence levels have the following relationships: 
$$\uppi_j[\boldsymbol{V}_n(j,j)]-\uppi_j(\bar{\upkappa}^{j,j})=\uppi_j^{(1)}(\bar{\upkappa}^{j,j})[\boldsymbol{V}_n(j,j)-\bar{\upkappa}^{j,j}]+o^2(1)$$
$$\uppi_j[\boldsymbol{\bar{F}}_n^{-1}(j,j)]-\uppi_j(\bar{\upkappa}^{j,j})=\uppi_j^{(1)}(\bar{\upkappa}^{j,j})[\boldsymbol{\bar{F}}_n^{-1}(j,j)-\bar{\upkappa}^{j,j}]+\tilde{O_d}^2(n^{-1/2})$$
$$\uppi_j[\boldsymbol{\bar{H}}_n^{-1}(j,j)]-\uppi_j(\bar{\upkappa}^{j,j})=\uppi_j^{(1)}(\bar{\upkappa}^{j,j})[\boldsymbol{\bar{H}}_n^{-1}(j,j)-\bar{\upkappa}^{j,j}]+\tilde{O_d}^2(n^{-1/2})$$
$$\uppi_j[\boldsymbol{\bar{H}}_n^{-1}(j,j)]-\uppi_j[\boldsymbol{\bar{F}}_n^{-1}(j,j)]=\uppi_j^{(1)}(\bar{\upkappa}^{j,j})(\boldsymbol{\bar{H}}_n^{-1}(j,j)-\boldsymbol{\bar{F}}_n^{-1}(j,j)+\tilde{O_d}^2(n^{-1/2})$$

$Proof$: First, we prove the differentiability of the confidence level function  $\uppi_j(x)$ (defined in (II.C)) to provide the necessary condition for the use of first and second order Taylor's expansion. Functions $f_1(x)=2x-1$, $f_{2j}(x)=C_jx^{1/2}$, and $\Phi$ are all infinitely differentiable on $(0,+\infty)$. Thus, the confidence level function $\uppi_j(x)=f_1\{\Phi [f_{2j}(x)]\}$ is also infinitely differentiable on $(0,+\infty)$.\par
Now we consider the confidence level gaps. From (1.1), we know $\bar{\upkappa}^{j,j}$ converges to $\boldsymbol{\bar{F}}(\boldsymbol{\uptheta}^*)^{-1}(j,j)$  as the sample size $n$ goes to infinity. By A3, $\boldsymbol{\bar{F}}(\boldsymbol{\uptheta}^*)^{-1}(j,j)$ exists and is positive. Thus for a sufficiently large $n$, $\bar{\upkappa}^{j,j}>0$. Thus, it is legitimate to calculate Taylor expansion of $\uppi_j(x)$ around $\bar{\upkappa}^{j,j}$.\par
A second-order Taylor expansion of  $\uppi_j(x)$ around $\bar{\upkappa}^{j,j}$ yields $$\uppi_j(x)-\uppi_j(\bar{\upkappa}^{j,j})=\uppi_j^{(1)}(\bar{\upkappa}^{j,j})(x-\bar{\upkappa}^{j,j})+1/2\uppi_j^{(2)}(a)(x-\bar{\upkappa}^{j,j})^2$$
where $\uppi_j^{(i)},i=1,2$ represents the $i^{th}$ order derivative of $\uppi_j$; and point $a$ lies between $x$ and $\bar{\upkappa}^{j,j}$. Plugging in $x=\boldsymbol{}{V}_n(j,j)$ and using Lemma 2 we have
\begin{align*}
    \uppi_j[\boldsymbol{V}_n(j,j)]-\uppi_j(\bar{\upkappa}^{j,j}) &= \uppi_j^{(1)}(\bar{\upkappa}^{j,j})(\boldsymbol{V}_n(j,j)-\bar{\upkappa}^{j,j})+1/2\uppi_j^{(2)}(a_{\mathbf{V}})(\boldsymbol{V}_n(j,j)-\bar{\upkappa}^{j,j})^2\\
     &= \uppi_j^{(1)}(\bar{\upkappa}^{j,j})(\boldsymbol{V}_n(j,j)-\bar{\upkappa}^{j,j})+o^2(1)
\end{align*}
where $a_\mathbf{V}$ lies between $\boldsymbol{V}_n(j,j)$ and $\bar{\upkappa}^{j,j}$. Plugging in $x=\boldsymbol{\bar{H}}_n^{-1}(j,j)$ or $x=\boldsymbol{\bar{F}}_n^{-1}(j,j)$, using Lemma 4 we get the result in Lemma 5 and the remainder term is $\tilde{O_d}^2(n^{-1/2})$.\par
Also, it follows that
$$\uppi_j[\boldsymbol{\bar{H}}_n^{-1}(j,j)]-\uppi_j[\boldsymbol{\bar{F}}_n^{-1}(j,j)]=\uppi_j^{(1)}(\bar{\upkappa}^{j,j})(\boldsymbol{\bar{H}}_n^{-1}(j,j)-\boldsymbol{\bar{F}}_n^{-1}(j,j)+\tilde{O_d}^2(n^{-1/2})$$
so we get results of Lemma 5.

\subsection{Main result}
In this section, we prove the main result in two steps: first, we measure the ratio of MSEs of the approximate confidence levels (defined in (\rom{2}.C)) for finite sample size $n$ ; second, we consider the infinite case ($n\rightarrow \infty$) and obtain the conclusion in the asymptotic sense.\par
$\boldsymbol{Theorem \quad 1}$. Under certain conditions \cite[pp. 25\textendash28]{cao2013relative} including the existence of limits of some sequences, interchanges of differentiation and integration for some functions, Weak Law of Large Numbers (WLLN) and Dominated Convergence Theorem (DCT), for any given $j$, the relative accuracy of the asymptotic tail probabilities associated with specified confidence levels for the observed and expected FIM satisfy:
\begin{equation}
    \liminf\limits_{n\rightarrow \infty}\frac{\textrm{MSE}_{H_{jj}}}{\textrm{MSE}_{F_{jj}}}\geq 1,j=1,...p
\end{equation}
If condition in equation (5) is also satisfied, then strict inequality holds in (6).\par

$Proof$: Let $D_n$ be the difference between squared errors of approximate confidence levels, then we have
$$D_{jn} \overset{\mathrm{\Delta}}{=} \{\uppi_j[\boldsymbol{\bar{H}}_n^{-1}(j,j)]-\uppi_j[\boldsymbol{V}_n(j,j)]\}^2-\{\uppi_j[\boldsymbol{\bar{F}}_n^{-1}(j,j)]-\uppi_j[\boldsymbol{V}_n(j,j)]\}^2$$
using Lemma 5, we could rewrite $D_{jn}$ as a summation of 4 parts:
$$D_{jn}=[\uppi_j^{(1)}(\bar{\upkappa}^{j,j})]^2(d_{jn1}+d_{jn2}+d_{jn3}+d_{jn4})$$ and we check the 4 parts separately.\par
Consider the first part. We have $$d_{jn1}=[\boldsymbol{\bar{H}}_n^{-1}(j,j)+\boldsymbol{\bar{F}}_n^{-1}(j,j)-2\boldsymbol{V}_n(j,j)][\boldsymbol{\bar{H}}_n^{-1}(j,j)-\boldsymbol{\bar{F}}_n^{-1}(j,j)]$$ taking the expectation of both sides and using Lemma 1\textendash5 and DCT, we get $$E(d_{jn1})=E(A_{njj}^2)+o(n^{-1})$$
\par
Also we have $$d_{jn2}=[\tilde{O_d}^2(n^{-{1}/{2}})+o^2(1)][\boldsymbol{\bar{H}}_n^{-1}(j,j)-\boldsymbol{\bar{F}}_n^{-1}(j,j)]$$
through a similar process we get $$E(d_{jn2})=o(n^{-1})$$.\par
Similarly we have
$$d_{jn3}=[\boldsymbol{\bar{H}}_n^{-1}(j,j)+\boldsymbol{\bar{F}}_n^{-1}(j,j)-2\boldsymbol{V}_n(j,j)]\cdot \tilde{O_d}^2(n^{-{1}/{2}})$$
and then $$E(d_{jn3})=o(n^{-1})$$
\par
Similarly we have
$$d_{jn4}=[\tilde{O_d}^2(n^{-{1}/{2}})+o^2(1)]\cdot \tilde{O_d}^2(n^{-{1}/{2}})$$ and $$E(d_{jn4})=o(n^{-1})$$
\par
Based on the results above, consider the expectation of $D_{jn}$:
\begin{equation}
    E(D_{jn})=[\uppi_j^{(1)}(\bar{\upkappa}^{j,j})]^2[E(A_{njj}^2)+o(n^{-1})]
\end{equation}

Then by condition A3 \cite{cao2013relative} which ensure the existence of the sequence and interchange of the limits, and continuity of the function $[\uppi_j^{(1)}(x)]^2$, we have 
\begin{equation}
\begin{split}
    \lim_{n\to\infty}[\uppi_j^{(1)}(\bar{\upkappa}^{j,j})]^2 &= [\uppi_j^{(1)}(\bar{\boldsymbol{F}}^{-1}(\boldsymbol{\uptheta}^*)(j,j))]^2\overset{\mathrm{\Delta}}{=}c_1
\end{split}
\end{equation}
According to Lemma 4,
\begin{equation}
    nE(A_{njj}^2)=nE[\tilde{O_d}(n^{-1/2})]^2=E[\tilde{O_d}(1)]^2=O_+(1)
\end{equation}
Taking the limit of both sides in (7), using (8) and (9), we have 
\begin{equation}
    \liminf_{n\rightarrow \infty}nE(D_{jn})=c_1\liminf_{n\rightarrow \infty}[nE(A_{njj}^2)+o(1)]=c_1\liminf_{n\rightarrow \infty}[O_+(1)+o(1)]\geq0
\end{equation}
Moreover, if the condition in equation (5) is also satisfied, then we have $$\liminf_{n\rightarrow \infty}nE(A_{njj}^2)>0$$ so we have
\begin{equation}
    \liminf_{n\rightarrow \infty}nE(D_{jn})
    =c_1\liminf_{n\rightarrow \infty}nE(A_{njj}^2)+c_1\liminf_{n\rightarrow \infty}o(1)>0.
\end{equation}
Under all necessary conditions \cite[A1\textendash A8]{cao2013relative}, we can write (6) as 
\begin{equation}
    \begin{split}
    \liminf\limits_{n\rightarrow \infty}\frac{\textrm{MSE}_{H_{jj}}}{\textrm{MSE}_{F_{jj}}}&=
   \liminf_{n\rightarrow \infty}\frac{E\{\uppi_j[\boldsymbol{\bar{H}}_n^{-1}(j,j)]-\uppi_j[\boldsymbol{V}_n(j,j)]\}^2}{E\{\uppi_j[\boldsymbol{\bar{F}}_n^{-1}(j,j)]-\uppi_j[\boldsymbol{V}_n(j,j)]\}^2}\\&=\liminf_{n\rightarrow \infty}\frac{nE(D_{jn})}{E\{\uppi_j[\boldsymbol{\bar{F}}_n^{-1}(j,j)]-\uppi_j[\boldsymbol{V}_n(j,j)]}+1\\ &\geq 1
    \end{split}
\end{equation}

In addition, if the condition in (5) is also satisfied, strict inequality in (12) holds according to (11).\par

Theorem 1 suggests that asymptotically the expected estimation error of the approximate confidence level with the expected FIM ($\textrm{MSE}_F$) never exceeds the one with the observed FIM ($\textrm{MSE}_H$), as long as the sample size $n$ is large enough.

\section{Numerical study}

We illustrate the theory above in three distinct estimation settings: Gaussian mixture distributions, signal-plus-noise with non-identically distributed noise, and parameter identification for state-space models. All of these settings arise regularly in practice and in the literature. Mixture problems are thoroughly reviewed in \cite{10.2307/2030064} and \cite{https://doi.org/10.2307/2981482}, along with an application in information theory\cite{53528} and bivariate quantile estimation\cite{doi:10.1080/01621459.1992.10475269}. The signal-plus-noise problem with non-identical noise distributions arises in practical problems where measurements are collected with varying quality of information across the sample. Some practical implications are discussed in \cite{doi:10.1080/03610928108828137} relative to the initial conditions in a Kalman filter model and \cite{SPALL1990297} in the context of outlier analysis. Another very common and important model is the state-space model, whose framework has been successfully applied in engineering, statistics, computer science and economics to solve a broad range of dynamical systems problems. Settings involving the calculation of uncertainty bounds related to MLEs of the parameters in state-space models are discussed in \cite{30995} and \cite{62251}. 

\subsection{Gaussian Mixture distribution}

Let $\boldsymbol{X}=[X_1,X_2,...X_n]^T$ be an i.i.d. sequence with probability density function:
$$f(x,\boldsymbol{\uptheta})=\frac{\uplambda}{\sqrt{2\uppi \upsigma^2}} \text{exp}\left(-\frac{(x-\upmu_1)^2}{2\upsigma^2}\right)+\frac{1-\uplambda}{\sqrt{2\uppi \upsigma_2^2}} \text{exp}\left(-\frac{(x-\upmu_2)^2}{2\upsigma^2}\right)$$
where $\boldsymbol{\uptheta}=[\uplambda,\upmu_1,\upmu_2]^T$ and $\upsigma$ is known. There is no closed form for MLE in this case. We use Newton’s method to achieve numerical approximation of $\hat{\boldsymbol{\uptheta}}_n$, and the covariance matrix of $\hat{\boldsymbol{\uptheta}}_n$ is approximated by the sample covariance of 1000 values of $\hat{\boldsymbol{\uptheta}}_n$ from 1000 independent realizations of data. The analytical form of the true FIM is not attainable, but the closed form of the Hessian matrix is computable since the log likelihood function is continuous differentialble. In this case we approximate the true FIM by numerical integration of Hessian matrix on the range of $X$.\par
In this study, we consider three cases when $\boldsymbol{\uptheta}^*=[0.5,0,4]^T, n=50$, and $\boldsymbol{\uptheta}^*=[0.5,0,2]^T, n=100$, and  $\boldsymbol{\uptheta}^*=[0.5,0,1]^T, n=100$ and for all cases $\upsigma=1$. For the second and third case, we use a bigger sample size $n$ to allow for adequate information to achieve reliable MLE when two individual Gaussian distributions have more overlapping area. We also show a typical outcome of $\bar{\boldsymbol{H}}_n(\hat{\boldsymbol{\uptheta}}_n)^{-1}$ and $\bar{\boldsymbol{F}}_n(\hat{\boldsymbol{\uptheta}}_n)^{-1}$. To define `typical', we consider the distance of $\bar{\boldsymbol{H}}_n(\hat{\boldsymbol{\uptheta}}_n)^{-1}$ or $\bar{\boldsymbol{F}}_n(\hat{\boldsymbol{\uptheta}}_n)^{-1}$ to the true covariance matrix, which is $\left\| \bar{\boldsymbol{H}}_n(\hat{\boldsymbol{\uptheta}}_n)^{-1}-n\text{cov}(\hat{\boldsymbol{\uptheta}}_n) \right\|$ or $\left\| \bar{\boldsymbol{F}}_n(\hat{\boldsymbol{\uptheta}}_n)^{-1}-n\text{cov}(\hat{\boldsymbol{\uptheta}}_n) \right\|$, and the typical value is the matrix with median distance to the true covariance matrix through all outcomes. In this case, we use the Frobenius norm.

\begin{table}[htbp]
\caption{Numerical Results for the Gaussian Mixture Example}
\begin{center}
 \begin{tabular}{|c|c|c|c|} 
 \hline
  True Parameter & $\boldsymbol{\uptheta}^*=[0.5,0,4]^T$ & $\boldsymbol{\uptheta}^*=[0.5,0,2]^T$ & $\boldsymbol{\uptheta}^*=[0.5,0,1]^T$\\ [0.5ex]
 \hline
 Sample size & $n=50$ & $n=100$ & $n=100$\\
 \hline
 $n\text{cov}(\hat{\boldsymbol{\uptheta}}_n)$ &  $\begin{pmatrix}
  0.26 & 0.14 & 0.10 \\ 
  0.14 & 2.50 & 0.44 \\
  0.10 & 0.44 & 2.55
  \end{pmatrix} $ & $\begin{pmatrix}
  1.42 & 2.86 & 2.49 \\ 
  2.86 & 10.82 & 5.52 \\
  2.49 & 5.52 & 7.93
  \end{pmatrix} $ & $\begin{pmatrix}
  6.39 & 8.19 & 8.81 \\ 
  8.19 & 18.65 & 8.87 \\
  8.81 & 8.87 & 20.89
  \end{pmatrix} $\\
\hline

%$\textrm{Typical}\atop \bar{\boldsymbol{F}}_n(\hat{\boldsymbol{\uptheta}}_n)^{-1}$
$\mycom{\textrm{Typical}}{\bar{\boldsymbol{H}}_n(\hat{\boldsymbol{\uptheta}}_n)^{-1}}$
  &  $\begin{pmatrix}
  0.26 & 0.07 & 0.08 \\ 
  0.07 & 1.97 & 0.25 \\
  0.08 & 0.25 & 2.79
  \end{pmatrix} $ & $\begin{pmatrix}
  0.96 & 1.56 & 1.45 \\ 
  1.56 & 5.85 & 2.86 \\
  1.45 & 2.86 & 5.19
  \end{pmatrix} $ & $\begin{pmatrix}
  14.94 & 11.97 & 23.19 \\ 
  11.97 & 12.09 & 16.65 \\
  23.19 & 16.65 & 44.96
  \end{pmatrix} $\\
  \hline
 $\mycom{\textrm{Typical}}{\bar{\boldsymbol{F}}_n(\hat{\boldsymbol{\uptheta}}_n)^{-1}}$  &  $\begin{pmatrix}
  0.26 & 0.06 & 0.07 \\ 
  0.06 & 2.09 & 0.26 \\
  0.07 & 0.26 & 2.51
  \end{pmatrix} $ & $\begin{pmatrix}
  0.87 & 1.73 & 1.45 \\ 
  1.73 & 7.75 & 3.02 \\
  1.45 & 3.02 & 4.09
  \end{pmatrix} $ & $\begin{pmatrix}
  7.61 & 20.01 & 7.08 \\ 
  20.01 & 62.52 & 17.81 \\
  7.08 & 17.81 & 8.39
  \end{pmatrix} $\\
  
  \hline
$\frac{\textrm{MSE}_{H_{jj}}}{\textrm{MSE}_{F_{jj}}}\textrm{ } \textrm{for} \textrm{ } \uplambda$& $1.28$ &$1.24$ &$1.11$\\
\hline
$\frac{\textrm{MSE}_{H_{jj}}}{\textrm{MSE}_{F_{jj}}}\textrm{ } \textrm{for} \textrm{ } \upmu_1$& $1.22$ &$1.14$ &$1.18$\\
\hline
$\frac{\textrm{MSE}_{H_{jj}}}{\textrm{MSE}_{F_{jj}}}\textrm{ } \textrm{for} \textrm{ } \upmu_2$& $1.22$ &$1.17$ &$1.20$\\
  \hline
\end{tabular}
\label{tab1}
\end{center}
\end{table}

Based on values of $\bar{\boldsymbol{H}}_n(\hat{\boldsymbol{\uptheta}}_n)^{-1}$ and $\bar{\boldsymbol{F}}_n(\hat{\boldsymbol{\uptheta}}_n)^{-1}$ we can calculate $\textrm{MSE}_{H_{jj}}$ and $\textrm{MSE}_{F_{jj}}$ for each parameter, leading to the ratios of errors in confidence levels shown in Table \rom{1}. Note that results in the table, with all ratios $\text{MSE}_{H_{jj}}/\text{MSE}_{F_{jj}}>1$, are consistent with the theoretical result.\par

\subsection{Signal-Plus-Noise problem}
Let $\boldsymbol{X}=[\boldsymbol{X}_1,\boldsymbol{X}_2,...\boldsymbol{X}_n]^T$ be an i.n.i.d. sequence, and
$$\boldsymbol{X}_i\sim N(\boldsymbol{\upmu},\boldsymbol{\Sigma}+\boldsymbol{Q}_i)$$
where $\boldsymbol{\upmu}$ is the common mean vector across observations, $\boldsymbol{\Sigma}$ is the common part of the covariance matrices and $\boldsymbol{Q}_i$ is the covariance matrix of noise for observation $i$. In practice, the $\boldsymbol{Q}_i$ are known and $\boldsymbol{\uptheta}$ contains unique elements in $\upmu$ and $\boldsymbol{\Sigma}$. There is no closed form for the MLE in this case. We use Newton’s method to achieve a good numerical approximation of $\hat{\boldsymbol{\uptheta}}_n$, which leads to $\bar{\boldsymbol{H}}_n(\hat{\boldsymbol{\uptheta}}_n)^{-1}$, and $\bar{\boldsymbol{F}}_n(\hat{\boldsymbol{\uptheta}}_n)^{-1}$. We also show a typical outcome of $\bar{\boldsymbol{H}}_n(\hat{\boldsymbol{\uptheta}}_n)^{-1}$ and $\bar{\boldsymbol{F}}_n(\hat{\boldsymbol{\uptheta}}_n)^{-1}$. The definition of `typical' is same as the numerical example above. We consider the cases of scalar data and 4-dimensional data below. \par
\subsubsection{1-d case}
First, we consider a 1-d case where 
$$X_i\sim N(\upmu,\upsigma^2+q_i)$$
and $q_i=0.1\times\textrm{mod}(i,10)$, $n=1000$, $\boldsymbol{\uptheta}=[\upmu,\upsigma^2]^T$ and the true value of the parameter is $\boldsymbol{\uptheta}^*=[10,10]^T$.

\begin{table}[htbp]
\caption{Covariance Matrices of the 1-d Signal-plus-noise Example}
\begin{center}
 \begin{tabular}{|c|c|} 
 \hline
  True Parameter & $\boldsymbol{\uptheta}^*=[10,10]^T$\\
 \hline
 Sample size & $n=1000$\\
 \hline
 $n\text{cov}(\hat{\boldsymbol{\uptheta}}_n)$ &  $\begin{pmatrix}
  13.43 & 1.83 \\ 
  1.83 & 426.24 
  \end{pmatrix} $\\ 
\hline
$\mycom{\textrm{Typical}}{\bar{\boldsymbol{H}}_n(\hat{\boldsymbol{\uptheta}}_n)^{-1}}$  &  $\begin{pmatrix}
  14.34 & -9.34\times 10^{-3} \\ 
  -9.34\times 10^{-3} & 396.34 
  \end{pmatrix} $\\
  \hline
 $\mycom{\textrm{Typical}}{\bar{\boldsymbol{F}}_n(\hat{\boldsymbol{\uptheta}}_n)^{-1}}$  &  $\begin{pmatrix}
  14.37 & 0 \\ 
  0 & 396.84 
  \end{pmatrix} $\\

  \hline
\end{tabular}
\label{tab2}
\end{center}
\end{table}

\begin{table}[htbp]
\caption{MSEs of the 1-d Signal-plus-noise example}
\begin{center}
\begin{tabular}{|c|c|c|c|}
\hline
& $\textrm{MSE}_{H_{jj}}$& $\textrm{MSE}_{F_{jj}}$ & $\frac{\textrm{MSE}_{H_{jj}}}{\textrm{MSE}_{F_{jj}}}$ \\
\hline
$\upmu$&$6.00\times10^{-5}$&$6.00\times10^{-5}$ & 1.00 \\
\hline
$\upsigma^2$&$1.08\times10^{-4}$&$1.06\times10^{-4}$ & 1.02 \\
\hline
\end{tabular}
\label{tab3}
\end{center}
\end{table}

As a result, for 1-d case the ratios of errors in confidence levels are shown in Table \rom{3}. These numerical results are consistent with the theoretical result. Specifically, in this case we have $\bar{\boldsymbol{H}}_n({\boldsymbol{\uptheta}}^*)^{-1}(1,1)=\bar{\boldsymbol{F}}_n({\boldsymbol{\uptheta}}^*)^{-1}(1,1)$ according to the analytical form of the log-likelihood function. \par

\subsubsection{4-d case}
Next we consider a 4-d case where 
$$\boldsymbol{X}_i\sim N(
\boldsymbol{\upmu},\boldsymbol{\Sigma}+\boldsymbol{Q}_i)$$
where $\boldsymbol{Q}_i=\sqrt{i}\times \boldsymbol{U}\boldsymbol{U}^T$ and
$\boldsymbol{U}$ is a $4\times4$ matrix where each entry is drawn from $Uniform(0,0.1)$ distribution. \par
Also in this case $\boldsymbol{\upmu}=[\upmu_1,\upmu_2,\upmu_3,\upmu_4]^T$ and $\boldsymbol{\Sigma}=\textrm{diag}\{ {\Sigma}_{11},{\Sigma}_{22},{\Sigma}_{33},{\Sigma}_{44} \}$, $\boldsymbol{\uptheta}=[\upmu_1,\upmu_2,\upmu_3,\upmu_4,{\Sigma}_{11},{\Sigma}_{22},{\Sigma}_{33},{\Sigma}_{44}]^T$ and the true value of the parameter is $\boldsymbol{\uptheta}^*=[0,0,0,0,1,1,1,1]^T$. Analogous to Table \rom{2} and Table \rom{3}, the ratios (all equal to or larger than 1) are shown in Table \rom{4}.

\begin{table}[htbp]
\caption{Numerical Results for the 4-d Signal-Plus-Noise Example}
\begin{center}
\begin{tabular}{|c|c|c|}
\hline

Sample size&$n=1000$ & $n=2000$  \\
\hline
$\frac{\textrm{MSE}_{H_{jj}}}{\textrm{MSE}_{F_{jj}}}\textrm{ } \textrm{for} \textrm{ } {\upmu}_{1}$&$1.00$ &$1.00$\\
\hline
$\frac{\textrm{MSE}_{H_{jj}}}{\textrm{MSE}_{F_{jj}}}\textrm{ } \textrm{for} \textrm{ } {\upmu}_{2}$&$1.00$ &$1.00$\\
\hline
$\frac{\textrm{MSE}_{H_{jj}}}{\textrm{MSE}_{F_{jj}}}\textrm{ } \textrm{for} \textrm{ } {\upmu}_{3}$&$1.00$ &$1.00$\\
\hline
$\frac{\textrm{MSE}_{H_{jj}}}{\textrm{MSE}_{F_{jj}}}\textrm{ } \textrm{for} \textrm{ } {\upmu}_{4}$&$1.00$ &$1.00$\\
\hline
$\frac{\textrm{MSE}_{H_{jj}}}{\textrm{MSE}_{F_{jj}}}\textrm{ } \textrm{for} \textrm{ } {\Sigma}_{11}$&$1.93$ &$3.26$\\
\hline
$\frac{\textrm{MSE}_{H_{jj}}}{\textrm{MSE}_{F_{jj}}}\textrm{ } \textrm{for} \textrm{ } {\Sigma}_{22}$&$2.22$ & $2.80$\\
\hline
$\frac{\textrm{MSE}_{H_{jj}}}{\textrm{MSE}_{F_{jj}}}\textrm{ } \textrm{for} \textrm{ } {\Sigma}_{33}$&$2.39$ & $1.55$\\
\hline
$\frac{\textrm{MSE}_{H_{jj}}}{\textrm{MSE}_{F_{jj}}}\textrm{ } \textrm{for} \textrm{ } {\Sigma}_{44}$&$1.83$ & $2.95$\\
\hline
\end{tabular}
\label{tab4}
\end{center}
\end{table}

Results in the above two examples are consistent with the theoretical conclusion in Section \rom{3}.B, i.e. asymptotically, for any element of the parameter, the confidence interval with the expected FIM is a better approximation than the one with the observed FIM in the MSE sense.

\subsection{State-space model}
A state-space model is a mathematical description of a physical system as a set of input, output and state variables. The state-space model considered is defined by the equations
$$\boldsymbol{x}_t=\boldsymbol{A}\boldsymbol{x}_{t-1}+\boldsymbol{w}_t$$
$$\boldsymbol{y}_t=\boldsymbol{C}\boldsymbol{x}_t+\boldsymbol{v}_t$$
for $t = 1,2,...,n$ time periods. And $\boldsymbol{x}_t$ is an unobserved $l$-dimensional state process, $\boldsymbol{A}$ is an $l\times l$ transition matrix, and $\boldsymbol{w}_t$ is a vector of $l$ zero-mean, independent disturbances with covariance matrix $\boldsymbol{Q}$. Let $\boldsymbol{y}_t$ be an observed $m$-dimensional process, $\boldsymbol{C}$ be an $m\times l$
design matrix, and $\boldsymbol{v}_t$ be a vector of $m$ zero-mean, independent disturbances with covariance matrix $\boldsymbol{R}$. The mean and covariance matrix of $\boldsymbol{x}_0$ (the initial $\boldsymbol{x}_t$) are denoted by $\boldsymbol{\upmu}$ and $\boldsymbol{P}_0$, respectively. It is assumed that $\boldsymbol{\upmu}$ and $\boldsymbol{\Sigma}$ are known and that $\boldsymbol{x}_0$, $\boldsymbol{w}_t$, and $\boldsymbol{v}_t$ are mutually independent and multivariate normal.\par
In our context, we consider situations where $\boldsymbol{A}$, $\boldsymbol{C}$, $\boldsymbol{R}$ are known. The unknown parameters of interest are the unique elements in the diagonal $\boldsymbol{Q}$. The calculation of log-likelihood function and Hessian matrix of the state-space model is discussed in detail in \cite{kitagawa2020computation}. In this example, we consider $l=3$, $m=1$, then we have the log-likelihood function of the system \cite{cao2013relative}:  
$$l(\boldsymbol{\uptheta})=-\frac{1}{2}\sum_{t=1}^n\textrm{log}(\boldsymbol{C}\boldsymbol{P}_{t|t-1}\boldsymbol{C}^T+{R})-\frac{1}{2}\sum_{t=1}^n(\boldsymbol{C}\boldsymbol{P}_{t|t-1}\boldsymbol{C}^T+{R})^{-1}\upepsilon_t^2$$
and the computation requires Kalman Filter equations:
$${\epsilon}_t={y}_t-\boldsymbol{C}\hat{\boldsymbol{x}}_{t|t-1}$$
$$\hat{\boldsymbol{x}}_{t|t-1}=\boldsymbol{A}\hat{\boldsymbol{x}}_{t-1|t-1}$$
$$\hat{\boldsymbol{x}}_{t|t}=\hat{\boldsymbol{x}}_{t|t-1}+\boldsymbol{K}_t\epsilon_t$$
$$\boldsymbol{K}_t=\boldsymbol{P}_{t|t-1}\boldsymbol{C}^T(\boldsymbol{C}\boldsymbol{P}_{t|t-1}\boldsymbol{C}^T+{R})^{-1}$$
$$\boldsymbol{P}_{t|t-1}=\boldsymbol{A}\boldsymbol{P}_{t-1|t-1}\boldsymbol{A}^T+\boldsymbol{Q}$$
$$\boldsymbol{P}_{t|t}=\boldsymbol{P}_{t|t-1}-\boldsymbol{K}_t\boldsymbol{C}\boldsymbol{P}_{t|t-1}$$
since in this case $R$, $y_t$ and $\epsilon_t$ are all scalar. With these equations in concert with the ``innovations form" of the Kalman filter \cite{spall93}, we can calculate the gradient and the Hessian matrix of the log-likelihood function \cite{kitagawa2020computation} in a recursive form. Also, the expected Fisher Information is attainable \cite{CAVANAUGH1996347}. \par
In this case, suppose $\boldsymbol{Q}=\text{diag}\{q_{11},q_{22},q_{33}\}$, $\boldsymbol{\uptheta}=[q_{11},q_{22},q_{33}]^T$ and specifically we have $\boldsymbol{A}=\begin{pmatrix}
0 & 1 & 0\\
0 & 0 & 1\\
0.8 & 0.8 & -0.8
\end{pmatrix}$, $\boldsymbol{C}=[1,0,0]$, ${R}=1$, $\boldsymbol{\upmu}=[0,0,0]^T$,  $\boldsymbol{P}_0=\begin{pmatrix}
0 & 0 & 0\\
0 & 0 & 0\\
0 & 0 & 0
\end{pmatrix}$, and the true value of parameter is $\boldsymbol{\uptheta}^*=[1,1,1]^T$. Since there is no closed form for the covariance matrix $n\text{cov}( \hat{\boldsymbol{\uptheta}}_n)$, in the following example, we use the sample covariance matrix\footnote{To get an idea of the reliability of the sample covariance matrix, we repeat the process of getting $n\text{cov}( \hat{\boldsymbol{\uptheta}}_n)$ for 100 times for the $n=100$ case. Then we get 100 sample covariance matrices, and we calculate the relative error for diagonal entries of these covariance matrices. To calculate the relative error, supppose the current sample covariance matrix of MLEs is $\boldsymbol{\Sigma}_0$, and the sample covariance matrix of another replication of MLEs is $\boldsymbol{\Sigma}_A$, then the component-wise relative error is $|\boldsymbol{\Sigma}_A(i,i)-\boldsymbol{\Sigma}_0(i,i)|/|\boldsymbol{\Sigma}_0(i,i)|,i=1,2,3$ for diagonal entries of the sample covariance matrix. As a result, the average relative error of each diagonal entry in $n\text{cov}( \hat{\boldsymbol{\uptheta}}_n)$ is 0.0723, 0.0893, 0.0860. It is desirable that we want a small relative error, which means the sample covariance matrix doesn't vary much when we have a different list of MLEs. Further theoretical and numerical methods are needed to get more accurate sample covariance matrices.} of 1000 $\hat{\boldsymbol{\uptheta}}_n$, and each $\hat{\boldsymbol{\uptheta}}_n$ is calculated from a sequence of observations by a stochastic search algorithm \cite[pp. 43\textendash 45]{spall_2003}. In this study we consider two cases: $\boldsymbol{\uptheta}^*=[1,1,1]^T, n=50$ and $\boldsymbol{\uptheta}^*=[1,1,1]^T, n=100$. We also show a typical outcome of $\bar{\boldsymbol{H}}_n(\hat{\boldsymbol{\uptheta}}_n)^{-1}$ and $\bar{\boldsymbol{F}}_n(\hat{\boldsymbol{\uptheta}}_n)^{-1}$, where the definition of `typical' is described in Section \rom{4}.A.

\begin{table}[ht]
\caption{Numerical Results for the State-Space Model}
\begin{center}
 \begin{tabular}{|c|c|c|} 
 \hline
  True Parameter & $\boldsymbol{\uptheta}^*=[1,1,1]^T$ & $\boldsymbol{\uptheta}^*=[1,1,1]^T$\\ [0.5ex]
 \hline
 Sample size & $n=50$ & $n=100$\\
 \hline
 $n\text{cov}(\hat{\boldsymbol{\uptheta}}_n)$ &  $\begin{pmatrix}
  0.3953 & -0.1598 & -0.1944 \\ 
  -0.1598 & 0.9062 & -0.1294 \\
  -0.1944 & -0.1294 & 0.7350
  \end{pmatrix} $  & $\begin{pmatrix}
  0.3342 & -0.1343 & -0.1957 \\ 
  -0.1343 & 0.3899 & -0.0269 \\
  -0.1957 & -0.0269 & 0.4470
  \end{pmatrix} $\\
\hline
$\mycom{\textrm{Typical}}{\bar{\boldsymbol{H}}_n(\hat{\boldsymbol{\uptheta}}_n)^{-1}}$  &  $\begin{pmatrix}
  1.9813 & -1.1514 & -1.0796 \\ 
  -1.1514 & 2.1761 & 0.2465 \\
  -1.0796 & 0.2465 & 1.0331
  \end{pmatrix} $ & $\begin{pmatrix}
  0.7762 & -0.4011 & -0.5309 \\ 
  -0.4011 & 0.6725 & 0.1249 \\
  -0.5309 & 0.1249 & 0.6204
  \end{pmatrix} $\\
  \hline
 $\mycom{\textrm{Typical}}{\bar{\boldsymbol{F}}_n(\hat{\boldsymbol{\uptheta}}_n)^{-1}}$  &  $\begin{pmatrix}
  1.0606 & -0.5137 & -0.7559 \\ 
  -0.5137 & 1.0009 & 0.0523 \\
  -0.7559 & 0.0523 & 1.1772
  \end{pmatrix} $ & $\begin{pmatrix}
  0.5493 & -0.2608 & -0.3965 \\ 
  -0.2608 & 0.4936 & 0.0356 \\
  -0.3965 & 0.0356 & 0.5949
  \end{pmatrix} $\\
  
  \hline
$\frac{\textrm{MSE}_{H_{jj}}}{\textrm{MSE}_{F_{jj}}}\textrm{ } \textrm{for} \textrm{ } q_{11}$& $1.0772$ &$1.2145$\\
\hline
$\frac{\textrm{MSE}_{H_{jj}}}{\textrm{MSE}_{F_{jj}}}\textrm{ } \textrm{for} \textrm{ } q_{22}$& $1.4841$ &$1.8824$\\
\hline
$\frac{\textrm{MSE}_{H_{jj}}}{\textrm{MSE}_{F_{jj}}}\textrm{ } \textrm{for} \textrm{ } q_{33}$& $1.5839$ &$1.6076$\\
  \hline
\end{tabular}
\label{tab5}
\end{center}
\end{table}
Results in Table \rom{5} show that for any element of the parameter, the approximate confidence level with the expected FIM is a better approximation than the one with the observed FIM in the MSE sense, and this is consistent with the theoretical result in Section \rom{3}.B.

\section{Conclusions and future work}
Confidence intervals are heavily used in the practice of parameter estimation in system identification. In this paper, we compare the confidence levels of the approximate confidence intervals for the components of the vector parameter with the observed and the expected FIM in the MSE sense. With theoretical proof and the support of numerical cases, it is suggested that the approximate confidence level corresponding to the expected FIM is higher than that corresponding to the observed FIM; thus, we conclude that under reasonable conditions, the approximate confidence intervals corresponding to the expected FIM outperforms the one corresponding to the observed FIM.\par

Although our research objective is a vector parameter, our current results are limited to properties of individual components of the vector. In future work, we will focus on the features of the confidence region of the entire vector parameter and sub-vector of the parameter vector. For example, under the asymptotic normality $\sqrt{n}(\boldsymbol{\hat{\uptheta}}_n-\boldsymbol{{\uptheta}}^*) \sim N(\boldsymbol{0},\boldsymbol{V}_n)$, $n(\boldsymbol{\hat{\uptheta}}_n-\boldsymbol{{\uptheta}}^*)\boldsymbol{V}_n^{-1}(\boldsymbol{\hat{\uptheta}}_n-\boldsymbol{{\uptheta}}^*)^{T}$ should asymptotically follow a chi-square distribution. So we could compare the performance of the observed FIM and the expected FIM replacing the true covariance matrix.\par
In the mixture Gaussian example, we could further explore the different behaviors of convergence between $\uplambda$, $\upmu_1$, $\upmu_2$. For example, when the true value of $\upmu_1$ and $\upmu_2$ are close to each other we need a larger sample size $n$ to get a persuasive numerical result. In the signal-plus-noise problem, we only discuss components in $\boldsymbol{\Sigma}$ as parameters because if components in $\boldsymbol{\uptheta}$ are all from the mean vector $\boldsymbol{\upmu}$, then we have the result $\boldsymbol{\bar{F}}(\boldsymbol{\uptheta}^*)^{-1}=\boldsymbol{\bar{H}}(\boldsymbol{\uptheta}^*)^{-1}$, and asymptotically $\textrm{MSE}_{H_{jj}}$ will be the same as $\textrm{MSE}_{F_{jj}}$. \par

This paper is based on two asymptotic assumptions, where the first is that the distribution of the MLE will converge to a normal distribution as the number of observations goes to infinity. The other one is that we use the inverse of the FIM based on samples and the MLE as an approximation for the real variance of MLE. Also, to estimate the covariance matrix of the MLE, we use the sample covariance matrix of the MLE as an approximation. The reliability of the sample covariance matrix of MLE needs to be further investigated, since it might vary with different samples. To make our results more accurate, we seek getting more stable sample covariance matrix, which doesn’t change much when we have different samples. In future work, we want to explore how well the sample covariance matrix could approximate the real covariance matrix both theoretically and numerically, which will solidify our theoretical foundation.

\end{document}